\documentclass[12pt]{amsart}
\usepackage{amscd,amsthm,amssymb,amsfonts,amsmath,epsfig,epic,eclbip,bbm,url}
\usepackage[arrow,matrix,tips,2cell]{xy}
\usepackage{color}

\theoremstyle{plain}
\newtheorem{theorem}{Theorem}[section]
\newtheorem{lemma}[theorem]{Lemma}
\newtheorem{proposition}[theorem]{Proposition}
\newtheorem{corollary}[theorem]{Corollary}

\theoremstyle{definition}
\newtheorem*{definition}{Definition}
\newtheorem{example}[theorem]{Example}
\theoremstyle{remark}
\newtheorem*{remark}{Remark}
\newtheorem*{acknowledgments}{Acknowledgments}

\newcommand{\Z}{\mathbb Z}    
\newcommand{\R}{\mathbb R}    
\newcommand{\C}{\mathbb C}    
\newcommand{\PP}{\mathbb P}   

\newcommand{\SL}{\operatorname{SL}}

\newcommand{\bsd}{\operatorname{bsd}}

\newcommand{\Star}{\operatorname{Star}}

\newcommand{\Conv}{\operatorname{Conv}}
\newcommand{\NC}{\operatorname{NC}}
\newcommand{\Lin}{\operatorname{Lin}}

\newcommand{\conv}{\operatorname{Conv}}
\newcommand{\relint}{\operatorname{relint}}
\newcommand{\aff}{\operatorname{aff}}

\newcommand{\suchthat}{\ : \ }
\newcommand{\<}{\langle}   
\renewcommand{\>}{\rangle} 

\newcommand{\smooth}{\Sigma \setminus D}
\newcommand{\D}{\Delta}
\newcommand{\Dv}{\D^\vee}

\newcommand{\E}{\nabla}
\newcommand{\Ev}{\nabla^\vee}
\newcommand{\am}{\mathcal A}

\begin{document}

\title[Integral affine structures III]{Integral affine structures on spheres 
\makebox[0pt]{\raisebox{.5in}{\normalfont DUKE-CGTP-05-03}}
III: complete intersections}

\author{Christian Haase}
\address{Mathematics Department \\ Duke University \\ Durham, NC 27708
  \\ USA}
\email{haase@math.duke.edu}
\author{Ilia Zharkov}
\address{Mathematics Department \\ Harvard University \\ Cambridge, MA 02138
  \\ USA}
\email{zharkov@math.harvard.edu}

\thanks{The first author was partially supported by NSF-grant
  DMS-0200740. The second author was partially supported by NSF grant
  DMS-0405939.} 

\begin{abstract}
  We extend our model for affine structures on toric Calabi-Yau hypersurfaces
  \cite{HZh} to the case of complete intersections.
\end{abstract}

\maketitle

\section{Introduction} \noindent
Starting  with the following data we will construct a pair of affine
structures on a sphere (or a product of spheres) with a codimension 2
discriminant locus and whose monodromy representations have dual
linear parts. Let $\D=\D^{(1)}+\dots+\D^{(r)}$ be a
nef-partition of a $d$-dimensional {\em reflexive} polytope $\D
\subset (\R^d)^*$. That is for each $i=1,\dots,r$ the polytope
$\D^{(i)}$ contains the origin and there is a {\em convex integral\/}
PL function $\psi_i$ which has values 1 on every nonzero vertex of
$\D^{(i)}$ and 0 on the other polytopes $\D^{(j)}$.

Let $\E=\E^{(1)}+\dots+\E^{(r)} \subset \R^d$ be the dual
nef-partition (also $d$-dimensional) as in \cite{BB}. Explicitly, if
$\phi_i$ is the Legendre dual to the zero function on $\D^{(i)}$,
i.e.,
$$\phi_i(x)=\max_{y\in\D^{(i)}}\{\<y,x\>\},$$ 
then the polytopes $\E^{(i)}$ are defined as:
$$\E^{(i)}:=\Conv \{0,\ x\in \Dv\ \suchthat \phi_i(x)=1\}.$$
\begin{figure}[htbp]
  \centering
  \includegraphics{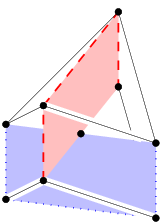}
  \qquad
  \includegraphics{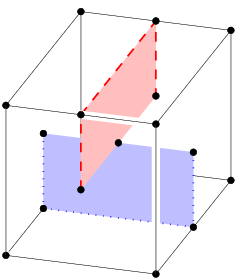}
  \hfill
  \includegraphics{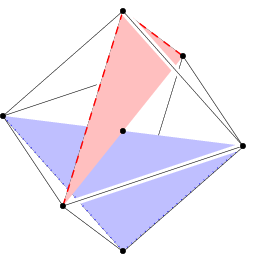}
  \qquad
  \includegraphics{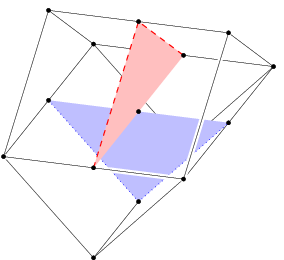}
  \\[2mm]
  $\E^\vee = \conv ({\color[rgb]{1,0,0}\D^{(1)}},
  {\color[rgb]{0,0,1}\D^{(2)}})$ \hfill
  $\D^\vee = \conv ({\color[rgb]{1,0,0}\E^{(1)}},
  {\color[rgb]{0,0,1}\E^{(2)}})$ \hfill \mbox{}
  \\
  \hfill $\D = {\color[rgb]{1,0,0}\D^{(1)}} +
  {\color[rgb]{0,0,1}\D^{(2)}}$
  \hfill \hfill $\E = {\color[rgb]{1,0,0}\E^{(1)}} +
  {\color[rgb]{0,0,1}\E^{(2)}}$
  \caption{Dual nef-partitions.}
  \label{fig:nefPartitions}
\end{figure}
(Here and later we will always denote by $P^\vee$ the polytope dual to
$P$.) In particular, the vertices of $\D^{(i)}$ are the gradients of
the piece-wise linear function $\phi_i$, and the vertices of
$\E^{(j)}$ are the gradients of the function $\psi_j$. Then we have
the duality (cf. \cite{Borisov}):
$$ \Ev=\Conv\{\D^{(1)},\dots,\D^{(r)}\} \text{ and } \Dv=\Conv\{\E^{(1)},\dots,\E^{(r)}\}.$$
The other part of our input are
two functions defined on the lattice points of $\Ev$ and $\Dv$:
$$\omega :\Ev_\Z \to \R, \qquad \nu: \Dv_\Z \to \R,$$
such that the induced subdivisions of $\Ev$ and $\Dv$ are central,
i.e. every maximal cell contains the origin. We will denote by $S$ and
$T$ the induced subdivisions of the boundaries of $\Ev$ and $\Dv$
respectively. In order to simplify the notation, we will add constants
to $\omega$ and $\nu$ so that $\omega(0)=\nu(0)=0$.

If $\omega$ and $\nu$ are integer valued then one can construct a
one-parameter family of complete intersections in a toric 
variety by taking the closure of the affine complete intersection
defined by $r$ Laurent polynomials in $(\C^*)^d$:
 $$f_i(z)=\sum_{m\in (\D^{(i)})_\Z} \lambda^{\omega(m)} z^m$$
in the toric variety $X_T$ defined by the fan over the subdivision
$T$. The function $\nu$ gives an integral K\"ahler class on $X_T$
which restricts to a class on the complete intersection. The
constructed affine structures are conjectured to be the ones which
arise in the metric collapse of the corresponding family.

We do not discuss the geometry of degenerations here (see a recent
article by Mark Gross~\cite{Gross} for more details on that aspect).
In the present paper we clarify the combinatorics of
the model. In particular, we establish a homeomorphism between the
model and a sphere. The affine structure and monodromy calculations
agree with those in \cite{Gross}.

\begin{acknowledgments}  
  We are indebted to Lev Borisov who provided several key ideas. We
  are also very grateful to Mark Gross for communicating to us Example
  \ref{schoen} which showed that the affine structure can be extended
  beyond the na\"{i}ve bipartite cover.
 
   Sphericity of our model was discussed with Richard Ehrenborg
  and Margaret Readdy, and with Vic Reiner. Frank Lutz, Niko Witte and
  G\"unter Ziegler came up with a counterexample to an overly
  optimistic conjecture of the first author.
\end{acknowledgments}

\section{The model}
\subsection{Semi-simplicity of nef-partitions}
A nef-partition $\D=\D^{(1)}+\dots+\D^{(r)}$ is called {\em reducible}
(cf. \cite{BB}) if  there exists a proper subset $\{i_1,\dots,i_s\}$
of the set $[r]=\{1,\dots,r\}$ such that the polytope
$\D^{(i_1)}+\dots+\D^{(i_s)}$ contains $0$ in its relative
interior. Otherwise the nef-partition is {\em irreducible}.
\begin{theorem}[\cite{BB}]
  Any nef-partition is a direct sum of irreducibles.
\end{theorem}

\begin{remark}
  In this decomposition one may need to refine the sum lattice so that
  it contains the direct sum of the constituent lattices as a finite
  index sublattice.
\end{remark}

From now on we will restrict our attention to irreducible
nef-partitions. Only those will correspond to true Calabi-Yau
families. Direct sums will correspond to products of Calabi-Yau
complete intersections of smaller dimensions, up to a finite group
action (cf. \cite{BB}).

\subsection{The sphere}
Consider the product $\D \times \E \subset (\R^d)^* \times \R^d$. The
complex $\Sigma$ --- our prospective sphere --- will be a
subdivision of
\begin{equation} \label{eq:sphere}
  |\Sigma| = \{ (m,n) \in \D \times \E \suchthat \< m,n \> = r \}.
\end{equation}
The rest of this subsection will be devoted to defining this subdivision.
Given $\sigma \in S$ we set $\sigma^{(i)} := \sigma \cap \D^{(i)}$,
and more generally $\sigma^I := \conv(\sigma^{(i)} \suchthat i \in
I)$ for a set $I \subseteq [r]$ of indices. Also denote $I_\sigma :=
\{ i \suchthat \sigma^{(i)} \neq \emptyset \}$.
We will say that $\sigma$ is {\em transversal} if $I_\sigma = [r]$,
that is, if $\sigma^{(i)}$ is not empty for every $i=1,\dots,r$. The
transversal cells form an upper order ideal $P$ in the face lattice of
$S$. For $\sigma \in P$ we will consider the subset $\sigma_\D$ of
$\D$:
$$\sigma_\D:=\sigma^{(1)}+\dots+\sigma^{(r)}.$$
(For non-transversal $\sigma$ this yields the empty set.)
The collection $\{\sigma_\D\suchthat \sigma\in P\}$ will be denoted
$S_\D$.

\begin{proposition} \label{prop:S_D}
  The collection $S_\D$ forms a polytopal complex whose face lattice
  is isomorphic to the poset $P$. In particular, the vertices of
  $S_\D$ are $\sigma_\D$ for minimal $\sigma\in P$.

  The space $|S_\D|\subset \partial\D$ coincides with the image of the
  natural projection map $p_1: |\Sigma| \to \D$.
\end{proposition}

\begin{proof}
  Here is an alternative definition of $S_\D$ from which the first
  assertion is obvious. The space $|S_\D|$ is the intersection of $\D$
  with the boundary of the dilation $r \Ev$. The cells of the complex
  $S_\D$ are $\sigma_\D = r\sigma \cap \D$.

  For the second assertion, if $x = p_1(x,y)$ with $(x,y) \in
  |\Sigma|$, then $\<\frac{1}{r}x,y\>=1$ shows that $x \in
  \partial(r\Ev)$. So $x \in |S_\D|$.
  Conversely, if $x \in \sigma_\D$, let $y \in \Dv$ with
  $\<\sigma,y\>=1$. Then $(x,y) \in |\Sigma|$.
\end{proof}

\begin{figure}[htbp]
  \centering
  \includegraphics{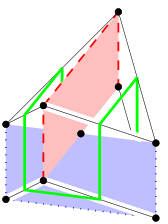}
  \qquad
  \includegraphics{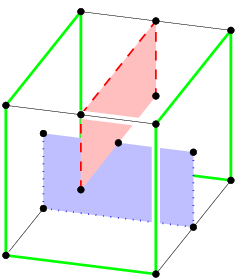}
  \hfill
  \includegraphics{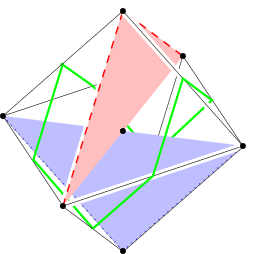}
  \qquad
  \includegraphics{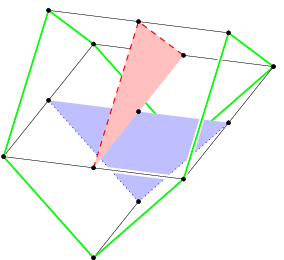}
  \caption{$|{\color[rgb]{0,1,0}S_\D}| = \partial(r\Ev) \cap \D$.}
  \label{fig:S_D}
\end{figure}

\begin{lemma}\label{lemma:sigma^i}
  For every (not necessarily transversal) $\sigma \in S$ and any index
  set $I \subseteq [r]$, $\sigma^I$ is a face of $\sigma$.
  If the two faces $\sigma^I$ and $\sigma^{\bar I}$ are non-empty,
  then the lattice distance between them is one.
  In particular, minimal transversal cells $\sigma$ are unimodular
  $(r-1)$-simplices.
\end{lemma}

\begin{proof}
  The integral PL function $\sum_{i \in I} \phi_i$ is linear on
  $\sigma$, equals 1 on $\sigma^I$ and vanishes on $\sigma^{\bar I}$.
  
  If $\sigma$ is a minimal transversal cell, then $\sigma_\D$ is a
  vertex so that each $\sigma^{(i)}$ is a vertex. Thus, $\sigma$ is a
  simplex which has facet width one by the first assertion.
\end{proof}

We can repeat everything said above for the dual data $(\E,\E^{(j)},T,
\tau^{(j)})$ instead of $(\D,\D^{(i)},S, \sigma^{(i)})$. In
particular, the poset $Q$ will be the collection of all transversal
cells $\tau \in T$. The subdivision $T_\E:=\{\tau_\E\suchthat \tau\in
Q\}$ of the relevant part of the boundary of $\E$ will coincide with
the image of the projection $p_2: |\Sigma| \to \E$.

We say that a pair $(\sigma,\tau)\in P\times Q$ is {\em adjoint} if
$\<\sigma_\D,\tau_\E\>=r$. Since $\<\D^{(i)},\E^{(j)}\>\le
\delta_{ij}$ a pair $(\sigma,\tau)$ is adjoint if and only if
$\<\sigma^{(i)},\tau^{(j)}\>=\delta_{ij}$.

\begin{definition}
  The complex $\Sigma$ consists of the cells $\sigma_\D \times
  \tau_\E$ for all adjoint pairs $(\sigma,\tau)$.
\end{definition}

\begin{remark}
  Everything said above is valid for reducible nef-partitions as
  well. If $\D$ is a direct sum $\D_1\oplus\dots\oplus\D_k$ of $k$
  irreducibles, then $\Sigma=\Sigma_1\times\dots\times\Sigma_k$.
  Later we will prove that for an irreducible nef-partition $\Sigma$
  is homeomorphic to the $(d-r)$-dimensional sphere; for reducible
  ones it is a product of spheres.
\end{remark}

\subsection{Tropical amoebas}
The function $\omega :\Ev_\Z \to \R$ and its restrictions \mbox{$\omega_i :
\D^{(i)}_\Z \to \R$} define tropical amoebas $\am$ and $\am^{(i)}$ as 
follows. Consider the $\omega$-lifted polytopes in $(\R^d)^* \times
\R$ 
\begin{equation*}
  \tilde \E^\vee = \conv \left[ \binom{m}{\omega(m)} \suchthat m \in
    \Ev_\Z \right] \text{ and } \tilde \D^{(i)} = \conv \left[
    \binom{m}{\omega(m)} \suchthat m \in \D^{(i)}_\Z \right] .
\end{equation*}
The lower convex hull of $\tilde \E^\vee$ projects to the subdivision
$S \star 0$ of $\Ev$. The normal fan of $\tilde \E^\vee$ subdivides
$\R^d \cong \R^d \times \{-1\}$ into cells
\begin{align*}
  F_\sigma \ &= \ \NC_{\tilde \E^\vee}(\tilde \sigma) \ \cap \ (\R^d
  \times \{-1\}) \\
  &= \ \{y \suchthat \<m,y\> - \omega(m) = \max_{m' \in \Ev_\Z}
  \<m',y\> - \omega(m') \text{ for all vertices } m \text{ of } \sigma
  \} .
\end{align*}
Similarly, for $\sigma \subseteq \D^{(i)}$,
\begin{align*}
  F^{(i)}_\sigma \ &= \ \NC_{\tilde \D^{(i)}}(\tilde \sigma) \ \cap \
  (\R^d \times \{-1\}) \\
  &= \ \{y \suchthat \<m,y\> - \omega(m) = \max_{m' \in \D^{(i)}_\Z}
  \<m',y\> - \omega(m') \text{ for all vertices } m \text{ of } \sigma
  \} . 
\end{align*}
The amoebas $\am$ and $\am^{(i)}$ are the polyhedral subcomplexes of
the cells of positive codimension:
\begin{align*}
  \am &:= \bigcup \ \{ F_\sigma \suchthat \sigma \in S \star 0, \dim
  \sigma > 0 \} \\
  \am^{(i)} &:= \bigcup \ \{ F^{(i)}_\sigma \suchthat \sigma \in
  \D^{(i)}, \sigma \in S \star 0, \dim \sigma > 0 \}
\end{align*}

We will be interested in the subcomplex
$\E_\omega$ of $\am$ which consists of all cells
$F_{\bar\sigma}$ for transversal $\sigma \in P$, where $\bar\sigma =
\sigma \star 0$. Observe that
with this definition, the barycentric subdivision $\bsd(\E_\omega)$ is
isomorphic to the order complex of $P$, i.e., the subcomplex of
$\bsd(S)$ induced by barycenters of transversal cells.

\begin{proposition}
  $\E_\omega$ is equal to the complex of bounded cells of $\bigcap
  \am^{(i)}$.
\end{proposition}

\begin{proof}
  We will show that the complexes of bounded cells of $\am$ and of
  $\bigcup \am^{(i)}$ agree with the boundary of the polytope
  $F_{0}=\bigcap F^{(i)}_0$.

  As the normal fan of the Minkowski sum $\tilde \D = \sum
  \tilde \D^{(i)}$ is the common refinement of normal fans of the
  $\tilde \D^{(i)}$, the subdivision of $\R^d \cong \R^d \times
  \{-1\}$ by the normal fan of $\tilde \D$ is the subdivision by
  $\bigcup \am^{(i)}$.

  The polytope $\tilde \E^\vee$ generates a cone $C$ with apex
  $\binom{0}{0}$.
  Its proper faces project to the cones over the cells of $S$. 
  The lower convex hull of each of the $\tilde \D^{(i)}$ is contained
  in the boundary of $C$.

  The Minkowski sum $\tilde \D = \sum \tilde \D^{(i)}$ generates the
  same cone: on one hand
  $\tilde \D \supseteq \tilde \E^\vee$, on the other hand
  the Minkowski sum of subsets of $C$ is contained in $C$.
  Moreover, the lower convex hull of $\tilde \D$ is contained in the
  boundary of $C$ as well: the sums $\sum \tilde \sigma^{(i)} \subseteq
  \partial C$ for $\sigma \in S \star 0$ belong to the lower hull, and
  their projections, the $\sum \sigma^{(i)}$, cover $\D$.
  
  The lower convex hull of $\tilde \D$
  projects to a mixed subdivision of $\D=\sum \D^{(i)}$. By the above
  observations, this is the subdivision of $\D$ by the cones over
  cells of $S$. That is, the cells are precisely the $\sum
  \sigma^{(i)}$ where $\sigma$ runs over the cells of $S \star 0$.

  \begin{figure}[htbp]
    \centering
    \includegraphics{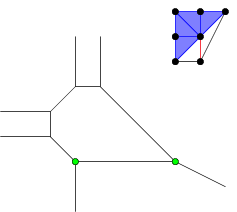}
    \qquad
    \includegraphics{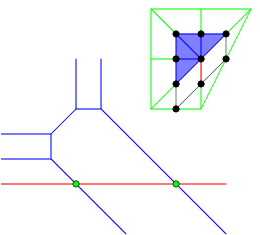}
    \caption{$\am$ versus $\bigcup \am^{(i)}$.}
    \label{fig:amoebas}
  \end{figure}

  Therefore, the bounded cells of both subdivisions are the
  intersection of the normal fan of $C$ with $\R^d \times \{-1\}$,
  and we have $F_{\bar\sigma} = \bigcap
  F^{(i)}_{\bar\sigma^{(i)}}$ for $\sigma \in S$.
\end{proof}

\begin{theorem} \label{thm:sphere}
  $\E_\omega$ is homeomorphic to the $(d-r)$-sphere.
\end{theorem}

We subdivide the proof into several lemmas which might be of
independent interest.

\begin{lemma} \label{lemma:v_i}
  There are vectors $v_i \in \relint \Delta^{(i)}$ so that
  $\sum v_i = 0$, and similarly $w_i \in \nabla^{(i)}$ so that $\sum
  w_i = 0$. The $v_i$ and the $w_i$ positively span $(r-1)$-dimensional
  linear spaces $V$ and $W$ respectively, and $(\R^d)^* = V \oplus
  W^\perp$. 
\end{lemma}

\begin{proof}
  We can write the origin as a convex combination (that is, a linear
  combination whose non-negative coefficients sum to $1$) of the
  vertices of $\Ev$ with all coefficients strictly positive. Define
  $v_i$ to be the contribution of the vertices of $\D^{(i)}$ to this
  sum.

  If we denote  $w^I := \sum_{i \in I} w_i$ for index sets $I
  \subseteq [r]$, then we see that $w^{\bar I}$ is non-positive on
  $\D^{(i)}$ for $i \in I$ so that $w^I = - w^{\bar I}$ is
  non-negative on $\D^{(i)}$ for $i \in I$.

  Because the nef-partition is irreducible, no proper subcollection
  of the $v_i$ or of the $w_i$ contains zero in their convex hull.
  Hence we have $(-1)^{\delta_{ij}} \<v_i,w_j\> < 0$, and any proper
  subcollection of the $v_i$ spans a pointed cone. The remaining
  assertions follow.
\end{proof}

Now shift the $\Delta^{(i)}$: set $\Delta^{(i)}_s = \Delta^{(i)} +
v_i$, and for $\sigma \in S$ set $\sigma_s = \conv(\sigma^{(i)}+v_i)$.

\begin{lemma} \label{lemma:shift}
  The $\sigma_s$ subdivide the boundary of $\nabla^\vee_s =
  \conv(\Delta^{(1)}_s, \ldots, \Delta^{(r)}_s)$ with the same
  combinatorics as $S$.
\end{lemma}

\begin{proof}
  We need to show that $\sigma_s$ belongs to the boundary of $\Ev_s$.
  The integral PL functions $\phi_j$ restrict to linear functions $y_j$
  on $\sigma$. We have $y_j(\sigma^{(i)})=\delta_{ij}$, and
  $y_j(\D^{(i)}) \le \delta_{ij}$. Thus $y=\sum y_j$ defines a face
  of $\Ev$ containing $\sigma$. 
  
  Furthermore, the $\lambda_{ij}:=\<y_j,v_i\>$ satisfy
  $(-1)^{\delta_{ij}} \lambda_{ij} < 0$, and $\sum_i \lambda_{ij} = 0$.
  Thus, the $r \times r$ matrix $\Lambda+\operatorname{id} =
  (\lambda_{ij} + \delta_{ij})_{1 \le i,j \le r}$ is invertible and we
  can find coefficients $\alpha_j$ so that for all $i$, $\alpha_i +
  \sum_j \lambda_{ij} \alpha_j = 1$. We claim that $\alpha_j \ge 0$.

  Assume $\alpha_j > 0$ for $j \in J_+$, $\alpha_j = 0$
  for $j \in J_0$, and $\alpha_j < 0$ for $j \in J_-$.
  We have 
  \begin{equation*}
    \sum_j \alpha_j = \mathbbm{1}^t \alpha =
    \mathbbm{1}^t(\Lambda+\operatorname{id}) \alpha = \mathbbm{1}^t
    \mathbbm{1} = r,
  \end{equation*}
  where $\mathbbm{1}$ is the all-one-vector. So $J_+ \neq \emptyset$.
  Consider the (row) vector $\beta$ with $\beta_j=|J_-|$ for $j \in
  J_+$, $\beta_j=0$ for $j \in J_0$, and $\beta_j=-|J_+|$ for $j \in
  J_-$. Then
  \begin{equation*}
    \beta (\Lambda+\operatorname{id}) \alpha = \sum_j \beta_j = 0 .
  \end{equation*}
  On the other hand, for $j \in J_+$, the $j^{\rm th}$ component of 
  $\beta (\Lambda+\operatorname{id})$ is
  \begin{equation*}
    |J_-| \underbrace{\left( \lambda_{jj} + \sum_{i
          \in J_+ \setminus j} \lambda_{ij} \right)}_{\ge 0} - |J_+|
    \underbrace{\sum_{i \in J_-} \lambda_{ij}}_{\le 0} + |J_-|> 0
    \text{ unless } J_- = \emptyset .
  \end{equation*}
  Similarly, for $j \in J_-$, the $j^{\rm th}$ component of 
  $\beta (\Lambda+\operatorname{id})$ is
  \begin{equation*}
    |J_-| \underbrace{\sum_{i \in J_+} \lambda_{ij}}_{\le 0} - |J_+|
    \underbrace{\left( \lambda_{jj} + \sum_{i \in J_- \setminus j}
        \lambda_{ij} \right)}_{\ge 0} - |J_+| < 0 
    \text{ as } J_+ \neq \emptyset .
  \end{equation*}
  From these equations we deduce the contradiction $\beta
  (\Lambda+\operatorname{id}) \alpha > 0$ unless $J_- = \emptyset$,
  that is $\alpha \ge 0$.

  We define $y_s := \sum \alpha_i y_i$ to get $\<y_s,\sigma^{(i)}_s\>
  = 1$, and $\<y_s,\D^{(i)}_s\> \le 1$.
\end{proof}

\begin{proof}[Proof of Theorem \ref{thm:sphere}]
  The face lattice of $\nabla_\omega$ is opposite isomorphic to the
  lattice $P$ of transversal cells of $S$.

  We claim that $W^\perp$ intersects precisely the (shifted) transversal
  cells so that the induced subdivision of the
  boundary of the convex $(d-r+1)$-polytope $\nabla^\vee_s \cap W^\perp$
  has face lattice $P$.
  
  \begin{figure}[htbp]
    \centering
    \includegraphics{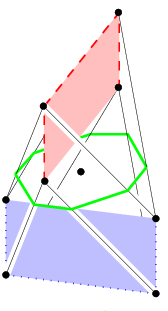}
    \qquad \qquad
    \includegraphics{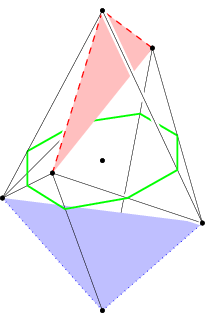}
    \caption{$\nabla^\vee_s \cap W^\perp$ and $\Dv_s \cap V^\perp$.}
    \label{fig:sphere}
  \end{figure}
  
  \begin{itemize}
  \item Suppose $\sigma \in S$ is not transversal, that is $I_\sigma
    \subsetneq [r]$.
    Then $\sigma_s$ can be separated from $W^\perp$ by the linear
    functional $w^{I_\sigma}$.
  \item Conversely, suppose $\sigma_s$ can be separated from $W^\perp$
    by some linear functional $w$. That is, $w \in (W^\perp)^\perp$
    and $w>0$ on $\sigma_s$. Then $w$ can be written as a positive
    combination of a proper subset of the $w_i$. But if $w_j$ is not
    in this collection, then $w<0$ on $\sigma^{(j)}_s$. Therefore
    $\sigma^{(j)}_s$ is empty, and $\sigma$ was not transversal.
  \end{itemize}
\end{proof}

\begin{corollary}
  If $\D=\D^{(1)}+\dots+\D^{(r)}$ is an irreducible nef-partition, then
  the topological space $\Sigma$ is homeomorphic to the sphere
  $S^{d-r}$.
\end{corollary}
\begin{proof}
  There is a (pulling) subdivision of $F_0$ that restricts
  to a subdivision of $\E_\omega$ isomorphic to the barycentric
  subdivision $\bsd(\Sigma)$. (Compare \cite[\S2.2]{HZh}.)
\end{proof}

\section{Affine structure on $\Sigma$}
\subsection{The smooth part and the discriminant}
The combinatorial structure of the affine (smooth) part of $\Sigma$ is
more complex in the complete intersection case than it is for 
hypersurfaces. The main difference is that the affine
structure can be extended beyond the bipartite covering. In
particular, the nerve of this enlarged covering is not a graph and the
fundamental group is not free. 

For every adjoint pair $(\sigma,\tau)$ we take the open star
neighborhood in $\bsd\Sigma$:
$$U_{(\sigma,\tau)}:=\Star(\hat\sigma_\D\times\hat\tau_\E).$$ 
The collection $\{U_{(\sigma,\tau)}\}$ provides an open 
covering of $\Sigma$.

For each minimal transversal $\sigma\in P$ it is convenient to combine
all charts 
$U_{(\sigma,\tau)}$ into a single chart
$U_\sigma$. Analogously, $V_\tau:=\bigcup_\sigma U_{(\sigma,\tau)}$ 
for a fixed minimal transversal
$\tau\in Q$. In terms of the two projections
\begin{equation*}
  p_1 \colon \bsd(\Sigma) \rightarrow \bsd(S_\D), \quad
  p_2 \colon \bsd(\Sigma) \rightarrow \bsd(T_\E),
\end{equation*}
$U_\sigma$ and $V_\tau$ are the preimages of the open star
neighborhoods of the vertices $\sigma_\D=\hat\sigma_\D$ and
$\tau_\E=\hat\tau_\E$ in the barycentric subdivision of $S_\D$ and
$T_\E$ respectively:
$$U_\sigma= p_1^{-1} (\Star(\sigma_\D)),\qquad V_\tau= p_2^{-1}
(\Star(\tau_\E)).$$
Then the closures of the $U_\sigma$'s cover the whole $|\Sigma|$ and the
combinatorial structure of this covering is dual to the poset
$P$. Similarly, $\{\bar V_\tau\}$ is a covering of $|\Sigma|$ dual to
the poset $Q$.

We say that that an adjoint pair $(\sigma,\tau)$ is {\em smooth} if
$$\dim \sigma^{(i)} \cdot \dim \tau^{(i)}=0, \qquad \text{ for all }
\ i=1,\dots,r.$$
Obviously, a pair $(\sigma,\tau)$ is smooth if either $\sigma$ or
$\tau$ is minimal. The {\em discriminant locus} $D$ is the full
subcomplex of $\bsd(\Sigma)$ generated by the vertices
$(\hat\sigma_\D\times\hat\tau_\E)$ for all {\em non-smooth} adjoint
pairs $(\sigma,\tau)$. Next we will describe an affine structure on
$\smooth$.

First, we fix a homeomorphism $\phi:\Sigma\to\E_\omega$, which sends
$U_\sigma$ to the corresponding (maximal) cell $F_{\bar\sigma}$ of
$\E_\omega$. (Recall that the $\bar\sigma = \sigma \star 0$ for
transversal $\sigma$ parameterize the cells of $\E_\omega$.)
For a smooth pair $(\sigma,\tau)$ we choose a partition
$\{I,J\}$ of the set $[r]$ such that $\dim \sigma^{(i)}=0$
for all $i\in I$ and $\dim \tau^{(i)}=0$ for $j\in J$.

Note that $\phi(U_{(\sigma,\tau)}) \subset \aff F_{\bar\sigma^I}
= \{x \in \R^d \suchthat \<m,x\>-\omega(m)=0,\ m \in
\sigma^I\}$. Then the integral affine structure on $U_{(\sigma,\tau)}$
is induced via $\phi$ by taking the quotient of a subspace in $\R^d$:
$$\aff F_{\bar\sigma^I}/\tau^J .$$
Any other allowable partition $\{I',J'\}$ will give a {\em
  canonically} equivalent structure. The easiest way to see this is to
choose the partition $\{I'',J''\}:=\{I\cap I',J\cup J'\}$:
\begin{align*}
  \aff F_{\bar\sigma^I}/\tau^J &\to \aff F_{\bar\sigma^{I''}}/\tau^{J''} \\
  \aff F_{\bar\sigma^{I'}}/\tau^{J'} &\to \aff F_{\bar\sigma^{I''}}/\tau^{J''}
\end{align*}
The condition $\<\sigma^{(i)},\tau^{(i)}\>=1$ ensures that the above
maps are isomorphisms of {\em integral\/} affine spaces.

In particular, the integral affine structure on $U_\sigma$ is induced
via the inclusion
$$\aff F_{\bar\sigma} \subset \R^d$$ 
from the standard integral affine structure on $\R^d$. The integral
affine structure on $V_\tau$ comes from the quotient $\R^d/\tau$. The
identification of the affine coordinates on the overlaps is provided
via the quotient map.
 
\subsection{Monodromy and extension of the affine structure}
The nerve of the partial covering $\{U_\sigma,V_\tau\}$ is the
bipartite graph $\Gamma$ whose nodes are labeled by the corresponding
minimal transversal cells $\sigma\in P$ and $\tau\in Q$.
We call a loop in $\Gamma$ {\it primary} if it has exactly 4 nodes
$\sigma_0, \tau_0, \sigma_1, \tau_1$
so that $\conv(\sigma_0,\sigma_1) \in S$ and $\conv(\tau_0,\tau_1) \in
T$. We denote such a loop by $(\sigma_0 \tau_0 \sigma_1 \tau_1)$, and
think of it as an element (possibly trivial) of the fundamental group
$\pi_1(\smooth)$ with a base point in $U_{\sigma_0}$. The
primary loops generate $\pi_1(\Gamma)$, and hence $\pi_1(\smooth)$, 
but there are many relations. 

\begin{proposition}
  The monodromy transformation 
  $T_{(\sigma_0 \tau_0 \sigma_1 \tau_1)} \colon \aff F_{\bar\sigma_0}
  \rightarrow \aff F_{\bar\sigma_0}$ 
  along the loop $(\sigma_0 \tau_0 \sigma_1 \tau_1)$ is given by
  $$T_{(\sigma_0 \tau_0 \sigma_1 \tau_1)} (x) = x + \sum_{j=1}^r 
  \left[\<\sigma^{(j)}_1, x\> - \omega(\sigma^{(j)}_1)\right]
  (\tau^{(j)}_1 - \tau^{(j)}_0).$$
\end{proposition}
\begin{proof}
  The vectors $\sigma^{(i)}_{0},\sigma^{(i)}_{1}\in
  (\R^d)^*$ and $\tau^{(j)}_{0},\tau^{(j)}_{1}\in \R^d$  satisfy
  \mbox{$\<\sigma^{(i)}_{0,1},\tau^{(j)}_{0,1}\>=\delta_{ij}$.} Hence,
  if $x\in \aff F_{\bar\sigma_0}$, then 
  $$x + \sum_{j=1}^r \left[\<\sigma^{(j)}_1, x\> -
    \omega(\sigma^{(j)}_1)\right] (\tau^{(j)}_1 - \tau^{(j)}_0) \in
  \aff F_{\bar\sigma_0} .$$
  Now put
  $$x' := x - \sum_{j=1}^r \left[\<\sigma^{(j)}_1, x\> -
    \omega(\sigma^{(j)}_1)\right] \tau^{(j)}_0 \in \aff
  F_{\hat\sigma_1} .$$
  Then 
  $x' \equiv x \mod \tau_0$ together with 
  $$x' \equiv  x + \sum_{j=1}^r \left[\<\sigma^{(j)}_1,x\> -
    \omega(\sigma^{(j)}_1)\right] (\tau^{(j)}_1 - \tau^{(j)}_0) \mod \tau_1$$ 
  imply the desired formula.
\end{proof}

We see that if $(\sigma,\tau)$ is a smooth pair then the monodromy
along any primary loop $(\sigma_0 \tau_0 \sigma_1 \tau_1)$ for minimal
transversal $\sigma_0, \sigma_1\subset\sigma$ and
$\tau_0,\tau_1\subset\tau$ is trivial. Thus, the affine structure
could be first defined on the $U_\sigma$'s and $V_\tau$'s and then
extended to the neighborhoods of the smooth vertices. This phenomenon
does not occur in the hypersurface case. We cannot, however, extend the
affine structure across the non-smooth $(\hat\sigma,\hat\tau)$ since
there are primary loops in the neighborhoods with non-trivial
monodromy. (Compare Corollary~\ref{cor:non-smooth}.)

\begin{corollary}
In a neighborhood of a vertex in $D$ the linear part of
the monodromy
\begin{equation*}
 \Lin(T) \colon \pi_1(\Star_\Sigma (\hat\sigma_\D\times\hat\tau_\E)
 \backslash D) \rightarrow \SL_{d-r}(\Z)
\end{equation*}
in a suitable basis is represented by a subgroup of the (abelian) group of
matrices in the form  
\begin{equation*}
  \begin{pmatrix}\operatorname{id}&*\\0&\operatorname{id}\end{pmatrix}
\end{equation*}
(the $(d-r)\times(d-r)$-identity matrix plus a $(\dim(\sigma)-r+1)
\times (\dim(\tau)-r+1)$ block in the upper right corner). 
\end{corollary} 
\begin{proof}
Fix a minimal transversal cell $\sigma_0 \subset \sigma$, and
consider $\Lin(T)$ as an endomorphism of the tangent space
$\sigma_0^\perp$ of $\aff F_{\bar\sigma_0}$.
The fundamental group $\pi_1(U_{(\sigma,\tau)}
 \setminus D)$ is generated (not freely) by the primary
loops $(\sigma_0 \tau_p \sigma_k \tau_q)$ for minimal $\sigma_k \subset  \sigma , \tau_{p,q} \subset  \tau$. 

Let $W_\tau \subset \sigma_0^\perp$ be the tangent space of $\tau_\D$,
i.e., the linear span (of dimension $\dim\tau-r+1$) of the
$\{\tau^{(j)}_p - \tau^{(j)}_q\}_{j=1}^r$, all minimal $\tau_{p,q} \subset
\tau$.
If we extend an integral basis for $W_\tau$ to an integral basis of
$\sigma_0^\perp$, then, because $W_\tau \subset \sigma_k^\perp$, the
linear part
\begin{equation*}
  \Lin(T_{(\sigma_0 \tau_p \sigma_k \tau_q)})(x) = x + \sum_{j=1}^r
  \<\sigma^{(j)}_k,x\> (\tau^{(j)}_p-\tau^{(j)}_q)
\end{equation*}
of the monodromy along any primary loop $(\sigma_0 \tau_p \sigma_k
\tau_q)$ will have the desired form. 
\end{proof}

\begin{corollary} \label{cor:non-smooth}
  For a primary loop $(\sigma_0 \tau_0 \sigma_1 \tau_1)$, the
  following are equivalent.
  \begin{enumerate}
  \item The monodromy around $(\sigma_0 \tau_0 \sigma_1 \tau_1)$ is
    trivial.
  \item There is a smooth pair $(\sigma,\tau)$ with $\sigma_0,
    \sigma_1 \subset \sigma$ and $\tau_0, \tau_1 \subset \tau$.
  \item For some $j$, $\sigma_0^{(j)} \neq \sigma_1^{(j)}$ and
    $\tau_0^{(j)} \neq \tau_1^{(j)}$.
  \end{enumerate}
\end{corollary}

\begin{proof}
  We only need to show that the monodromy around non-smooth loops is
  non-trivial.
  Because $\sigma_0$ and $\sigma_1$ belong to a common face of $S$, we
  can find an $x \in \sigma_0^\perp$ so that $\<\sigma_1, x\> \ge 0$,
  and $\<\sigma_1^{(j)}, x\> > 0$. Analogously, choose $y \in
  \tau_0^\perp$ so that $\<y, \tau_1\> \ge 0$, and $\<y,
  \tau_1^{(j)}\> > 0$.
  Then
  \begin{equation*}
    \<y, \Lin(T_{(\sigma_0 \tau_0 \sigma_1 \tau_1)})(x) - x\> > 0 ,
  \end{equation*}
  so that $\Lin(T_{(\sigma_0 \tau_0 \sigma_1 \tau_1)})(x) \neq x$.
\end{proof}

\begin{example}\label{schoen} (Compare \cite{Gross} for details).
Take the complete intersection in $\PP^1\times\PP^2\times\PP^2$ of two
hypersurfaces of tridegrees (1,3,0) and (1,0,3). In $\R^5$ and
$(\R^5)^*$ we let
\begin{equation*}
\begin{split}
\D_1=\Conv\{&(0,-1,-1,0,0),(0,2,-1,0,0),(0,-1,2,0,0),\\
          &(1,-1,-1,0,0),(1,2,-1,0,0),(1,-1,2,0,0)\},\\   
\D_2=\Conv\{&(0,0,0,-1,-1),(0,0,0,2,-1),(0,0,0,-1,2),\\
          &(-1,0,0,-1,-1),(-1,0,0,2,-1),(-1,0,0,-1,2)\},\\
\E_1=\Conv\{&(1,0,0,0,0),(0,-1,0,0,0),(0,0,-1,0,0),(0,1,1,0,0)\},\\
\E_2=\Conv\{&(-1,0,0,0,0),(0,0,0,-1,0),(0,0,0,0,-1),(0,0,0,1,1)\}.          
\end{split}
\end{equation*}
For simplicity we choose $\omega=\mathbbm{1}$ and
$\nu=\mathbbm{1}$. Then $\Sigma$ is $S^3$ and the discriminant locus
is the Hopf link of two circles each with multiplicity 12. The
complement $\smooth$ is homotopy equivalent to a 2-torus. The global
monodromy $\Z^2\to\SL_3(\Z)$ can be represented by abelian matrices:
\begin{equation*}
(a,b) \mapsto \left(
\begin{array}{ccc}
 1 & 12a & 12b\\
 0 & 1 & 0 \\
 0 & 0 & 1
\end{array}
\right)
\end{equation*} 
\end{example} 

\subsection{A glimpse of mirror symmetry}
The dual affine structure arises from interchanging the r\^oles of
$\D,\omega$ and $\E,\nu$. The topological space $\Sigma$ and the
discriminant remain literally the same. The linear part of the
monodromy is given by the transpose inverse matrices if one chooses
the dual bases. This constitutes the topological part of the duality.

The geometric part is much more involved. Recall that the affine
structure (not just its monodromy representation) depends on the map
$\phi:\Sigma\to\E_\omega$. We can also consider its dual counterpart -
the map $\phi^\vee:\Sigma\to\D_\nu$. To define an integral {\em
  K\"ahler} affine structure we need the composition $\phi^\vee \circ
\phi^{-1}: \E_\omega \to \D_\nu$ be locally potential. The affine
Calabi conjecture says that for every pair of polarizations
$(\omega,\nu)$ there is a unique choice for $\phi^\vee \circ
\phi^{-1}$ such that the local potentials solve the real
Monge-Amp\`ere equation.  

\bibliographystyle{alpha}
\bibliography{references}

\end{document}